\documentclass[11pt,reqno]{amsart}
\makeindex

\usepackage{amsmath,amsfonts,amssymb,amsthm,amscd}

\usepackage{amssymb}
\usepackage[latin1]{inputenc}
\usepackage{epsfig,epsf}
\usepackage{enumerate}
\usepackage{indentfirst}
\usepackage{euscript}

 \theoremstyle{plain}
 \newtheorem{theorem}{Theorem}
\newtheorem{proposition}{Proposition}
\newtheorem{lemma}{Lemma}

 \newtheorem{remark}{Remark}
 \theoremstyle{definition}

\title[Ricci  Tensors on ${\mathbb R}^n$]
{
Ricci  Tensors  with Rotational Symmetry on ${\mathbb R}^n$ }

\author{Ronaldo A. Garcia  and Romildo S. Pina}

 \thanks{ The first author is fellow of CNPq. 
 This work was done under the project PRONEX/FINEP/MCT - Conv. 76.97.1080.00 -
Teoria Qualitativa das Equa\c c\~oes Diferenciais Ordin\'arias and
CNPq, Grant 476886/2001-5}

  \thanks{{\sc key-words:} Ricci tensor, Rotational symmetry,  Potential function, Implicit differential equation.
  \;  Math. Sub. Classif. MSC: 53C21, 58J603, 4A09  }

\begin{document}
\maketitle
 \vskip .5cm

\begin{abstract}  In this paper   is considered  the differential equation $Ric(g)=T$,
where $Ric(g)$ is the Ricci tensor of the metric $g$ and $T$ is a   rotational
 symmetric tensor on ${\mathbb R}^n$.   A  new, geometric,  proof of the existence
    of smooth  solutions of this equation, based on qualitative theory of implicit
    differential equations,  is presented here. This result was obtained previously by
     DeTurck and   Cao in 1994.

\end{abstract}


\section{ Introduction  } \label{sec:introducao}
In this paper   will be considered a particular case of the  second order partial differential equation
\begin{equation} Ric(g) =T, \label{eq:ricci} \end{equation}

\noindent where $g$ is a Riemannian metric  in a manifold $\mathbb M^n$, $Ric(g)$ is the Ricci tensor of $g$
and $T$ is a given symmetric tensor of second order. This equation is
 of physical significance in {\it field theory}, see chapter XI of \cite{ll} and chapter 18 of   \cite{ta}.
 For example, the Einsten's gravitational equations, the Maxwell's equations of electromagnetic fields
 and  the Euler equations for fluids are
related to equation 1. The tensor $T$ is interpreted physically as the stress-energy tensor due to
 the presence of matter.

 DeTurck showed that if $M$ is a surface then the equation $Ric(g)=T$ can be solved
 locally if $T=\rho \gamma$ where $\rho$ is a smooth real function and $\gamma$ is a
  positive definite tensor, see  \cite{de1}.

Also   DeTurck showed  that in dimension 3 or more the problem $Ric(g)=T$, with $T$
 non singular, has local   solution in the smooth or analytic category and that for
  $T=x_1dx_1^2+dx_2^2+\cdots + dx_n^2$, singular at $x_1=0$, then the problem has no
  local solution near $x_1=0$, see \cite{de2}.

Recently, DeTurck and Goldschmidt  studied the equation \ref{eq:ricci} with $T$ singular,
 but of constant rank. A detailed analysis of  integrability conditions  for local
 solvability was obtained. Also the authors obtained various  results of existence of local
   solutions, under additional hypothesis, see  \cite{deturck}.

In this work    will be considered the   Ricci equation $Ric(g)=T_S$,
 where $g=e^{2f}g_0$ is a conformal deformation of the   canonical metric
 $g_0= dr^2+ r^2d\Theta^2$  of $\mathbb R^n$ and $T_S$ is a given tensor
 with   $SO(n)$ rotational symmetry.
This problem was considered previously by  DeTurck and  Cao, \cite{cade}.
 They obtained local solutions near the origin $0$.

 This paper
provides a new proof of   existence and unicity, up to homothety, of smooth
 local solutions near the origin of $\mathbb R^n$ for the equation $Ric(g)=T_S$.

In the paper by DeTurck and Cao, \cite{cade},  no explicit statement about
 the smoothness of the metric $g$ at $0$ is presented.

This paper is organized as follows. In section \ref{sec:pre}   the definition of Ricci tensor  is
 recalled and the problem $Ric(g)=T$ is formulated. In section \ref{sec:rst} the main results
  of this work are stated in proved. In section \ref{sec:hrs} we give an example of a metric $g$ with
  rotational symmetry in $\mathbb R^n$ and the associated Ricci tensor is calculated explicitly. Finally
  in section \ref{sec:cr} some general problems are stated.

\section{ Preliminaries } \label{sec:pre}

On a  n-dimensional Riemannian manifold $(M,g)$ the associated
 Riemann or {\it curvature tensor} $R=R(g)$   is given, in a local chart,
  by the coefficients $R_{ijkl}$.
The following relations hold.
$$\aligned R_{ijkl} &= R_{klij}\\
R_{ijkl}+R_{iljk}+ R_{iklj}&=0\\
\nabla_m R_{ijkl} +\nabla_k R_{ijlm}+\nabla_l R_{ijmk}&=0\endaligned$$

\noindent This space of coefficients  has dimension $n^2(n^2-1)/12$.

The  contraction $R_{ik}=g^{jl}R_{ijkl}$ of the  {\em curvature tensor}
 is the called the {\it Ricci tensor}, see \cite{aubin} and \cite{jost}.
 So, for each $p\in M$ the Ricci tensor is a symmetric bilinear form
  $Ric(g): T_pM\times T_pM\to \mathbb R$ and therefore the dimension
  of the space of coefficients $R_{ik}$ is equal to $n(n+1)/6$.
 The {\it Ricci curvature}
in the direction $X=\{X_i\}$ is $R_{ij}X_iX_j$ and $g^{ij}R_{ij} $ is
 called the {\it scalar curvature.}

In local coordinates the Ricci tensor is given by:

\begin{equation}
\aligned
R_{ij}=&\frac 1{2(n-1)} g^{kl}[\frac{\partial^2 g_{jl}}{\partial x^i\partial x^k}
+\frac{\partial^2 g_{ik}}{\partial x^j\partial x^l} -
\frac{\partial^2 g_{kl}}{\partial x^i\partial x^j}-\frac{\partial^2 g_{ij}}{\partial x^k\partial x^l}] \\
+&\frac 1{n-1} g^{kl}g_{pq}[\Gamma_{ik}^p\Gamma_{jl}^q -\Gamma_{ij}^p\Gamma_{kl}^q] \cr
  =&\frac{\partial \Gamma_{ij}^s}{\partial x^s}  -
  \frac{\partial \Gamma_{is}^s}{\partial x^j} +
  \Gamma_{ij}^s\Gamma_{st}^t -\Gamma_{it}^s\Gamma_{sj}^t\endaligned
\end{equation}

where, \begin{equation} \Gamma_{ij}^l =
 \frac 12 g^{kl}[\frac{\partial g_{jk}}{\partial x^i}
  +\frac{\partial g_{ ik}}{\partial x^j} -\frac{\partial g_{ij }}{\partial x^k}]
 \end{equation}

\noindent are  the Christoffel symbols of the metric $g$.

In the two dimensional situation $R_{1212}=\mathcal K$ is the
Gaussian curvature and it  is the only non zero coefficient of both tensors.

In the three dimensional case there exists an algebraic  relation
 between $R(g)$ and $Ric(g)$ which is given by:

$$R_{ijkl}=g_{ik}R_{jl}- g_{il}R_{jk}-
g_{jk}R_{il}+
g_{jl}R_{ik}-\frac 12 R( g_{ik}g_{jl}-g_{il}g_{jk}).$$
\noindent  Here $R=g^{ik}R_{ik}$ is the scalar curvature, see \cite{ha} and \cite{po}.

For $n\geq 4$, in general, the Ricci tensor does not determine the Riemannian curvature tensor.

Also the following holds.

Let $M\subset \mathbb R^{n+1}$ be a hypersurface with second fundamental form $h$
and $\{X_1,\cdots, X_n\}$ be an orthonormal frame given by the principal directions
 $X_i$. It follows  that the principal curvatures are given by $h_i=h(X_i,X_i)$ and
  according to \cite{po}   the Riemman tensor is given by :

$$\aligned R_{ijkl}=&h(X_i,X_k)h(X_j,X_l)-h(X_j,X_k)h(X_i,X_l)\\
=& h_ih_j( \delta_{ik}\delta_{jl}-\delta_{jk}\delta_{il}).\endaligned$$

So the Ricci tensor is given by
$$R_{ij}= \delta_{ij}[h_i(h_1+\cdots+h_n)-h_i^2].$$

The scalar curvature is equal to
$$R=(\sum_{i=1}^n h_i)^2- \sum_{i=1}^n h_i^2.$$

In this paper we will consider   the following restricted  problem, studied
   by  Cao and  DeTurck, \cite{cade}.

\noindent{\bf Problem:}{\em  Given a smooth rotational symmetric   tensor $T$
 on ${\mathbb R}^n$, determine, if it exists,  a smooth metric $g$ such that   }
\begin{equation} \label{eq:re} Ric(g)=T.\end{equation}

\noindent

\section{Rotationally  Symmetric Tensors}\label{sec:rst}

Consider a  tensor $T$ on the n-dimensional Euclidean space $\mathbb R^n$
  symmetric with respect to the orthogonal group $SO(n)$; that is $\gamma*T=T$ for
  every $\gamma\in SO(n)$. These tensors will be referred to as {\it  rotationally symmetric. }

Under the hypothesis of non singularity of $T$,  the following lemma was proved in \cite{cade}.

\begin{lemma} Let  $T=\varphi(t)dt^2+t^2\psi(t)d\Theta^2$, $t\in \mathbb R_+$ and
$\Theta\in {\mathbb S}^{n-1}$, a smooth, nonsingular and rotationally symmetric
tensor on $\mathbb R^n$. Then $T$ is either positive or negative definite everywhere.
Moreover  $\varphi(0)=lim_{t\to 0} \psi(t)$.
\end{lemma}

Consider also a smooth rotationally symmetric metric $g,$ expressed  in spherical
 coordinates $(r,\Theta)=(r,\theta_1,\cdots, \theta_{n-1})$, as

\begin{equation}\label{eq:gr}
g=2e^{f(t)}[r^\prime(t)dt^2+r(t)^2 d\Theta^2].
\end{equation}

\noindent where, $r(0)=0$, $r^\prime(0)=1$ and $r^\prime(t)>0$ for all $t>0$.

The Ricci tensor of $g$ defined  by equation \ref{eq:gr}, see \cite{cade}, is given by:

\begin{equation} \label{eq:rn} Ric(g)=\alpha(r)dr^2+ r^2\beta(r)d\Theta^2\end{equation}

\noindent where,

 $$\aligned \alpha(r) &=-(n-1)[f_{rr}+\frac{f_r}r]\\
 \beta(r) &=-[f_{rr}+(2n-3)\frac{f_r}r+(n-2)(f_r)^2]\\
 f_r &=f^\prime(t)/r^\prime(t)\;\;\;\;\text{ and}\;\;\;\;  f_{rr}=f^\prime_r(t)/r^\prime(t).\endaligned$$

Therefore the equation $Ric(g)=T$ is equivalent to the following.

 $$ t^2\psi=[\frac{\varphi}{n-1}(\frac{r}{r^\prime})^2]-(n-2)[2rf_r+(rf_r)^2].$$

If $\varphi \ne 0$ it follows that:

\begin{equation}\label{eq:rt}
\frac{t^2\psi\varphi}{n-1}=[\frac{\varphi  r}{(n-1)r^\prime}]^2
-\frac{(n-2)}{n-1}\varphi[2rf_r+(rf_r)^2]\end{equation}

In \cite{cade}  the {\em Ricci potential} was defined  as

$$w(t)=\frac{1}{2\pi}\int_{D_t}KdA_g.$$

\noindent  Here $D_t=\{(s,\Theta)\in W : 0\leq s\leq t\}$ and $W$ is a fixed
two dimensional subspace of $\mathbb R^n$ and $K$ is the sectional curvature of
$W$. The function $w$ is well defined since the metric $g$ is rotationally symmetric.

\begin{lemma} \label{lem:rps} If $w$ is the Ricci potential of the rotationally
 symmetric metric $g=2e^{f(t)}[r^\prime(t)dt^2+r(t)^2 d\Theta^2]$ then
\begin{equation} \label{eq:rp}\aligned
w(t)=&-rf_r=-\frac{r(t)}{r^\prime(t)}f^\prime(t)\\
w^\prime(t)=& \frac{\varphi r}{(n-1)r^\prime} \endaligned
\end{equation}
\end{lemma}

\begin{proof} The sectional curvature of the plane $W$ is given by $$K= -e^{-2f}[f_{rr}+\frac{f_r}r].$$
 Therefore it follows that
$$ w(t)= \frac{1}{2\pi}\int_{D_t}KdA_g= \frac{1}{2\pi}\int_0^t \int_0^{2\pi} K(t) r(t)e^{2f(t)}r^\prime (t)dt.$$
Then,

$$\aligned \frac{dw}{dt}=& K(t) r(t)e^{f(t)}r^\prime (t)\\
=& - [f_{rr}+\frac{f_r}r]r(t) r^\prime (t) \\
=&-\frac{d}{dt}[rf_r].\endaligned $$

As $w(0)=0$ and $r(0)=0$ it follows that $w(t)=-r(t)f_r(t)$. Differentiating the equation
 above it follows that
$$w^\prime(t)=-r^\prime f_r-r f_{rr}r^\prime=\frac{\varphi(t)r(t)}{(n-1)r^\prime(t)}.$$
\end{proof}

In terms of the Ricci potential  $w$
 it follows, from equation \ref{eq:rt},   the following implicit differential equation.

\begin{equation}
\label{eq:iw} \left(\frac{dw}{dt}\right)^2=\frac{1}{n-1}[(n-2)\varphi(t)(w^2-2w)+t^2\varphi(t)\psi(t)]
\end{equation}

\begin{proposition} \label{prop:siw} The implicit differential equation
 \ref{eq:iw}
 with $\varphi(0)=\psi(0)\ne 0$, has a unique local smooth solution, with initial condition $w(0)=w^\prime(0)=0$,
defined in the  interval $[0,\epsilon)$, such that $w^\prime(t)\varphi(0)>0$ for all $t\in (0,\epsilon)$.
\end{proposition}

\begin{proof} For $n=2$ the implicit differential equation \ref{eq:iw} is  equivalent to the ordinary
 differential equations $dw/dt=\pm t \sqrt{\varphi(t)\psi(t)}$.
  Direct integration leads to the result
  stated.

So we suppose $n>2$.  Consider   smooth extensions of $\varphi$  and $\psi$ for $t<0$ and  the
implicit surface
$${\mathcal F}(t,w,p)=\frac{1}{n-1}[(n-2)\varphi(t)(w^2-2w)+t^2\varphi(t)\psi(t)]-p^2 \,\,=0, $$
\noindent where $p=\frac{dw}{dt}$.

Under the hypothesis above, we have $d{\mathcal F}(0)=[0 \;\;
 \frac{n-2}{n-1}\varphi(0) \;\; 0]\ne 0$, so ${\mathcal F}^{-1}(0)$ is
  locally a  regular smooth surface near $0$.

Next we consider the smooth Lie-Cartan vector field
$$X={\mathcal F}_p\frac{\partial}{\partial t}+p{\mathcal F}_p  \frac{\partial}{\partial w}
-({\mathcal F}_t+p{\mathcal F}_w)\frac{\partial}{\partial p},$$

\noindent defined in a tubular neighborhood of the surface ${\mathcal F}^{-1}(0)$ and tangent to it.

\noindent The projections of the integral curves of $X$ by $\pi(t,w,p)=(t,w)$ are
the solutions of  equation \ref{eq:iw}.

The origin $0$ is a singular point of $X$, isolated in ${\mathcal F}^{-1}(0),$ and

$$ \aligned DX(0)=\begin{pmatrix}0 & 0 & -2\\
0&0&0\\
\frac{-2\varphi(0)\psi(0)}{n-1}& -2\frac{(n-2)\varphi^\prime(0)}{n-1}
& 2\frac{(n-2)\varphi(0)}{n-1}
\end{pmatrix}\endaligned $$

The non zero eigenvalues of $DX(0)$, $\lambda_1$ and $\lambda_2$  are the roots of
 $$ -\lambda^2+ 2\frac{(n-2)\varphi(0)}{n-1}\lambda +\frac{4\varphi(0)\psi(0)}{n-1}=0.$$

As $\lambda_1\lambda_2= -4\frac{\varphi(0)\psi(0)}{n-1} <0$,  it follows that
$0$ is a hyperbolic saddle point of $X|_{{\mathcal F}^{-1}(0)}$. Also the    unstable
 and stable separatrices are transversal to the regular fold curve $\Sigma$, defined
 by $\Sigma= \{(t,w,p) | {\mathcal F}(t,w,p)={\mathcal F}_p(t,w,p)=0\}$.

Therefore the projections, by $\pi$,  of the smooth  stable and unstable separatrices
  of the hyperbolic singularity $0$  are smooth and quadratically tangent  to $\pi(\Sigma)$.
  Direct calculation shows that $\pi(\Sigma)=\{(t,w(t),0): w(t)=
 \frac{\psi(0)}{n-1}\frac{t^2}{2}+\cdots\}$. The projections of the separatrices of
 the hyperbolic saddle will be called {\it folded separatrices}. The integral
 curves of $X$ and of their projections in the plane $ (t,w)  $ are as shown
 in the Figure 1
below.
\end{proof}

\begin{figure}[htbp]\label{figure:fig1}
\begin{center}
\includegraphics[angle=0, width=12cm]{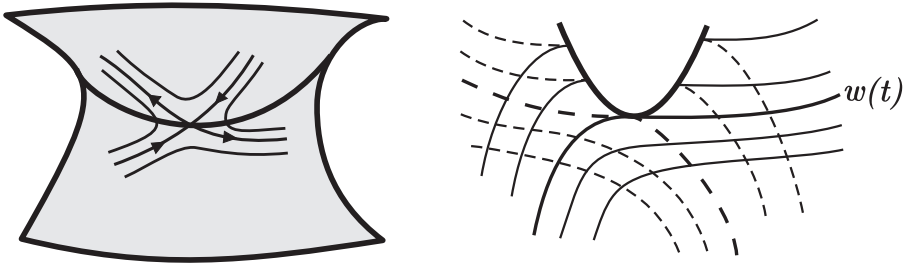}
       \caption { Folded saddle point }
   \end{center}
 \end{figure}

 \begin{remark} The proposition \ref{prop:siw} corresponds to Lemma 2.1
  (also labelled Proposition 2.1) and Propositions 2.2 and 2.3 of \cite{cade}.
   In the mentioned paper there is no statement about the smoothness of  the
    solution $w $ at $0$, the fixed point of the rotation group.
We adopt to work in the $C^\infty$ category, but with the obvious changes
  the result is valid in the $C^r$, $r\geq 2$, category.

\end{remark}

 \begin{remark} The analysis above is  similar to that carried out in the study
of asymptotic lines near a   parabolic curve, see \cite{gas}.  \end{remark}

\begin{lemma}\label{lem:pcs} Let $w$ be the smooth solution of equation (\ref{eq:iw}).
Consider the  singular differential equation

\begin{equation}\label{eq:pcs}\aligned
 (n-1)& w^\prime r^\prime -\varphi r=0\\
r(0)= &0, \;\; r^\prime(0)=1\\
w(0)= &\; w^\prime(0)=0, \; \; w^{\prime\prime}(0)=\frac{\varphi(0)}{n-1}  \endaligned
\end{equation}

Then there exists a smooth solution $r=r(t) $ of equation (\ref{eq:pcs}) in a interval $[0,\epsilon)$.

\end{lemma}

\begin{proof} Let $r(t)=t R(t)$ and $w(t)=\frac{t^2}2 W(t)$. Then it follows that the
 equation (\ref{eq:pcs}) is equivalent to the following equation

\begin{equation}\label{eq:rpc}
t\frac{dR}{dt}= \frac{\varphi R}{(n-1)( W+\frac t2 W^\prime)} -R
\end{equation}

Therefore, the line $(0, R)$ is a normally hyperbolic set and, by Invariant Manifold
 Theory, \cite{hps}, there exists a smooth solution $R(t)$ defined in a neighborhood
 of $0$ with initial condition $R(0)=1$.
\end{proof}

\begin{lemma}\label{lem:pcsf} Let $w$ be the smooth solution of equation (\ref{eq:iw}).
Consider the  singular differential equation
\begin{equation}\label{eq:pcsf}\aligned
 (n-1)& w^\prime f^\prime +w\varphi  =0\\
f(0)= &0,  \\
w(0)= &\; w^\prime(0)=0, \; \; w^{\prime\prime}(0)=\frac{\varphi(0)}{n-1}  \endaligned
\end{equation}

Then there exists a smooth solution $f=f(t)$ of equation (\ref{eq:pcsf}) in a interval $[0,\epsilon)$.

\end{lemma}

\begin{proof} The same argument as in the proof of lemma \ref{lem:pcs} works here.
\end{proof}

>From   proposition \ref{prop:siw} and lemmas \ref{lem:pcs} and \ref{lem:pcsf}
  follows the next proposition.

\begin{proposition} Let $T=\varphi(t)dt^2+t^2\psi(t)d\Theta^2$ be non singular
everywhere and suppose $w$ is a solution of equation (\ref{eq:iw}) such
that $w(0)=w^\prime(0)$ and $w^\prime(t)\varphi(t)>0 $ for $t>0$. Then
the Ricci system $Ric(g)=T$ is solvable. In fact, $ g=2e^{f(t)}[r^\prime(t)dt^2+r(t)^2 d\Theta^2]$,
 where  $r$ and $f$ are as stated, respectively, in lemmas \ref{lem:pcs} and \ref{lem:pcsf}.

Also formally we can write,

\begin{equation}
\aligned
 r(t) =& t  \text{ exp}\left( \int_0^t \left[ \frac{\varphi(s)}{(n-1)w^\prime(s)}-\frac 1s \right]ds\right),\\
f(t) =& -\int_0^t w(s)\frac{r^\prime(s)}{r(s)}ds +c=-\int_0^t \frac{\varphi}{n-1}\frac{w}{w^\prime}ds+c.\endaligned
\end{equation}

\end{proposition}

\begin{proof}
>From lemma \ref{lem:rps}   it follows that $w(t)=-r(t)f_r(t)=-r(t)f^\prime(t)/r^\prime(t)$
and $w^\prime(t)=\varphi(t) r(t)/(n-1)r^\prime(t)$.     From lemmas \ref{lem:pcs} and \ref{lem:pcsf}
   follows that    $r(t)$ and $f(t)$ are smooth solutions of these equations.
Therefore, integration  leads to the following result.
\end{proof}

\begin{theorem} Consider  the smooth, nonsingular, rotationally symmetric tensor
 $T=\varphi(t)dt^2+t^2\psi(t)d\Theta^2$. Suppose that ${\mathcal F}^{-1}(0)  $
 is a regular surface for  all $ t
\geq 0$ and $\frac{d}{dt}(t^2\psi(t))\varphi(t)\ne 0$, i.e., the set $\Sigma$,
 defined by $\Sigma= \{(t,w,p) | {\mathcal F}(t,w,p)={\mathcal F}_p(t,w,p)=0\}$ is a regular curve.
 Then $Ric(g)=T$ has a rotationally symmetric solution $g$ defined on all $\mathbb R^n$.
\end{theorem}

\begin{proof} The solution of the Ricci equation $Ric(g)=T$ is obtained from the stable
 or unstable separatrix of a hyperbolic saddle of $X$. This separatrix is  defined until
  it reaches the boundary of a connected component of the set
   $\{(t,w,p) |\;  {\mathcal F}(t,w,p)\geq 0\}$, which is, under the hypothesis above,
    the regular curve $\pi(\Sigma)$.
    The condition $\frac{d}{dt}(t^2\psi(t))\varphi(t)\ne 0$ means that the folded curve
     $\Sigma=\{(t,w,p): {\mathcal F}(t,w,p)={\mathcal F}_p(t,w,p)=0\}$ is a regular curve,
      with two connected components and that    the vector field $X$ has no singular point
       outside $0$ on the connected component of $\pi(\Sigma)$ that contains $0$. Therefore,
        the folded separatrices of the  saddle point $0$ of $X$ cannot reach the boundary
         of $\{(t,w,p) | \; {\mathcal F}(t,w,p)\geq 0\}$. If this  occurs there
would be a topological disk, bounded by a folded separatrix and by a connected component
 of the folded curve, foliated by regular curves transversal, outside $0$, to the folded
  curve. But this is impossible.
\end{proof}


\section{Hypersurfaces with Rotational Symmetry}\label{sec:hrs}

In this section we will calculated the Ricci tensor for a rotationally symmetric hypersurface
 of $\mathbb R^{n+1}$.

Let $\alpha: \mathbb R^n\to \mathbb R^{n+1}$ be an embedding with rotational symmetry, i.e,
a graph of a function $h$,  given by $\alpha(y_1,\cdots, y_n)= (y_1,\cdots, y_n, h(y_1^2+\cdots + y_n^2))$.

In  spherical coordinates it follows that:

$\alpha(r,\theta_1,\cdots, \theta_{n-1})=(y_1,\cdots, y_n,y_{n+1})$ where,

\begin{equation}\label{eq:gsi} \aligned y_1=& r\cos\theta_1\cdots \cos\theta_{n-1}\\
y_2=& r\cos\theta_1\cdots \cos \theta_{n-2}\sin\theta_{n-1}\\
\cdots =& \cdots \\
y_{n-1}=& r\cos\theta_{1}\sin\theta_2\\
y_n=&r \sin \theta_1\\
y_{n+1}=& h(r^2)\endaligned \end{equation}

Therefore the first fundamental form of $\alpha$ is given by $g=(g_{ij})$, where
$$\aligned g_{11}=&1+4r^2( h^\prime(r^2))^2\\
g_{22}=&r^2\\
g_{33}=& r^2\cos^2\theta_1\\
\cdots =&\cdots\\
g_{nn} =& r^2\cos^2\theta_1\cos^2\theta_2\cdots \cos^2\theta_{n-2}\\
g_{ij}=&0,\; i\ne j\endaligned $$

In a concise form we can write

$$g=[1+4r^2( h^\prime(r^2))^2]dr^2+r^2d\Theta^2,$$
\noindent where $d\Theta^2$ is the metric of the unitary sphere $\mathbb S^{n-1}$.

In the  diagonal metric $(g_{ij})$ above the Ricci tensor is given by

$$\aligned Ric(g)=&Ric(\frac{\partial}{\partial r},\frac{\partial}{\partial r})dr^2+
\sum_{i=1}^{n-1} Ric(\frac{\partial}{\partial \theta_i}, \frac{\partial}{\partial \theta_i})d\theta_i^2\\
 +& \sum_{i=1}^{n-1} Ric(\frac{\partial}{\partial r},\frac{\partial}{\partial \theta_i})drd\theta_i+
\sum_{i,j, i\ne j}^{n-1} Ric(\frac{\partial}{\partial \theta_i},\frac{\partial}{\partial \theta_j})d\theta_id\theta_j.
\endaligned$$

A long, but straightforward, calculation   gives:
$$\aligned
Ric(\frac{\partial}{\partial r},\frac{\partial}{\partial r})=&\frac{(n-1)f^\prime(r)}{2rf(r)}\\
Ric(\frac{\partial}{\partial \theta_1}, \frac{\partial}{\partial \theta_1})=&
\frac{rf^\prime(r)}{2f(r)^2}-\frac{n+2}{f(r)}+n-2   \\
Ric(\frac{\partial}{\partial \theta_i}, \frac{\partial}{\partial \theta_i})=&
 [\frac{rf^\prime(r)}{2f(r)^2}-\frac{n+2}{f(r)}+n-2]\prod_{k=1}^{i-1}\cos^2\theta_k,\; 2\leq i\leq n-1\\
Ric(\frac{\partial}{\partial \theta_i}, \frac{\partial}{\partial \theta_j})=&0\\
Ric(\frac{\partial}{\partial r}, \frac{\partial}{\partial \theta_i})=&0,
\endaligned $$
\noindent
where $f(r)=1+4r^2( h^\prime(r^2))^2$.

So    the following proposition holds.

\begin{proposition} Let $\alpha: \mathbb R^{n}\to \mathbb R^{n+1}$ be an embedding
with rotational $SO(n)$ symmetry , which in spherical coordinates is expressed   by
equation \ref{eq:gsi}. Then the Ricci tensor of the induced metric $g=(g_{ij})$, is given by:

$$\aligned Ric(g )=&a(r)dr^2+b(r)dr^2\\
=& \frac{(n-1)f^\prime(r)}{2rf(r)}dr^2+[\frac{rf^\prime(r)}{2f(r)^2}-\frac{n+2}{f(r)}+n-2]d\Theta^2.\endaligned$$
Here, $f(r)=1+4r^2( h^\prime(r^2))^2$.\end{proposition}

Finally we  remark that  the principal curvatures of the embedding
$\alpha$ are given by:

$$h_1=\frac{ 2h^\prime(r^2)+4r^2 h^{\prime\prime}(r^2)}{[1+4r^2( h^\prime(r^2))^2]^{3/2}},
\hskip 1cm h_2=\cdots =h_n=\frac{2 h^\prime(r^2)}{[1+4r^2( h^\prime(r^2))^2]^{1/2}}.$$

\section{Concluding Remarks}\label{sec:cr}

There is a considerable literature about the equation $Ric(g)=T$,
see \cite{de1}, \cite{de2}, \cite{deturck} and \cite{pt},  and the
general problem is the following.

\noindent{\bf Problem: }{\em  Given a    tensor $T$  on a Riemannian
 manifold $\mathbb M^n$, determine, if it exists,  a metric $g$ such that}

\begin{equation} \label{eq:reg} Ric(g)=T.\end{equation}

This equation is a second order system of quasi linear  partial differential equation,    \cite{hor}.

Other problems related to   the equation  $Ric(g)=T$ are  the following classical   Nirenberg   and Yamabe problems.

For $n=2$ consider the two-sphere $\mathbb S^2$ with the standard metric $g_0=dx^2+dy^2+dz^2$.

The Gaussian curvature of $g=e^{2u}g_0$ is given by
\begin{equation}\label{eq:1}
K(p)=(1-\Delta)e^{-2u(p)},
\end{equation}
\noindent where $\Delta$ is the Laplacian relative to the metric $g_0$.

A  global  problem in this case is the following: {\em which functions $K$ can be the
 Gaussian curvature of a metric $g$ which is a conformal deformation of $g_0$, i. e.,
 for which $K:{\mathbb S}^2\to \mathbb R$ are there solutions $u$ of equation (\ref{eq:1})?}

A   general version of this problem in $\mathbb R^n$, $n\geq 3$,  is known as the generalized
{\it Yamabe Problem} and consists in  obtaining  solutions of the partial differential equation

\begin{equation}\label{eq:2}
4{ \frac{n-1}{n-2}} \Delta_g u+ R_g u=R_{\bar g} u^{ \frac{n+2}{n-2}};  \;\;\; u>0,
\end{equation}
 where $\bar{g}= u^{4/(n-2)} g$, $R_g$ is the scalar curvature of $g$ and $R_{\bar g} $
 is the prescribed scalar curvature of the metric $\bar{g}$, see
 \cite{aubin}.

Another kind o problem is the local realization problem for the Gaussian curvature of a
 surface which can be stated as follows:  {\it  given a germ $K$ of a smooth function of
  two variables near the origin, find a surface in $\mathbb R^3$ with Gaussian curvature
  equal o $K$.} This problem was considered by Arnold, \cite{ar}, and the
  main result is that it can be solved whenever $K$ has a critical point of finite
  multiplicity at the origin.

Some  more concrete problems can be also stated.

\noindent{\bf Problem 1:} Existence and unicity of solutions for the equation  $ Ric(g)=T $
 in manifolds with boundary, for example in the unitary disk $\mathbb D^n\subset \mathbb R^n$
  or in  the cylinder $\mathbb D^n\times \mathbb R^m$.

\noindent{\bf Problem 2:} Study of the equation $Ric(g)=T$ in $\mathbb R^{m+n}$ where $T$
 has the symmetry of other geometric groups, for example $O(m)\times O(n)$. See \cite{cade}.

 \noindent{ \bf Problem 3:} In the singular case, i. e.,   $T=\varphi(t)dt^2+t^2\psi(t)d\Theta^2$,
 with $\varphi(0)=0$ and $\varphi^\prime(0)\ne 0$ analyze the existence and unicity of
  local solutions of the symmetric Ricci problem.

\noindent {\bf Problem 4:} Consider the Ricci principal curvatures  defined by
the equation $R_{ij}-\lambda g_{ij}=0$ and  the associated Ricci principal
directions. Study the {\em Ricci Configuration}, defined by  $n$ one dimensional
 singular foliations on a Riemannian manifold $(\mathbb M, g)$ and compare it
  with the {\it principal configuration} of a hypersurface  of $\mathbb R^{n+1}$.
  This setting is   analogous to that of the   configurations  of principal
  curvature lines, see \cite{gs} and \cite{garcia}.

 \newpage
\author{\noindent  Ronaldo Alves Garcia \\Instituto de Matem\'{a}tica e Estat\'{\i}stica,
\\Universidade Federal de Goi\'as,\\CEP 74001-970, Caixa Postal 131,\\Goi\^ania, GO, Brazil\\
\\
Romildo da Silva Pina \\Instituto de Matem\'{a}tica e Estat\'{\i}stica,
\\Universidade Federal de Goi\'as,\\CEP 74001-970, Caixa Postal 131,\\Goi\^ania, GO, Brazil
}
\end{document}